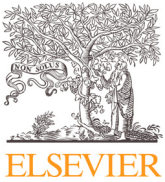

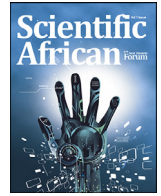

# The Sod gasdynamics problem as a tool for benchmarking face flux construction in the finite volume method

## Osama A. Marzouk

*College of Engineering, University of Buraimi, P.O. Box 890, Post Code 512 Al Buraimi, Sultanate of Oman*

### A B S T R A C T

The Finite Volume Method in Computational Fluid Dynamics to numerically model a fluid flow problem involves the process of formulating the numerical flux at the faces of the control volume. This process is important in deciding the resolution of the numerical solution, thus its quality. In the current work, the performance of different flux construction methods when solving the one-dimensional Euler equations for an inviscid flow is analyzed through a test problem in the literature having an exact (analytical) solution, which is the Sod problem.

The work considered twenty two flux methods, which are: exact Riemann solver (Godunov), Roe, Kurganov-Noelle-Petrova, Kurganov-Tadmor, Steger-Warming Flux Vector Splitting, van Leer Flux Vector Splitting, AUSM, AUSM$^+$, AUSM$^+$$-$up, AUFS, five variants of the Harten-Lax-van Leer (HLL) family, and their corresponding five variants of the Harten-Lax-van Leer-Contact (HLLC) family, Lax-Friedrichs (Lax), and Rusanov.

The methods of exact Riemann solver and van Leer showed excellent performance. The Riemann exact method took the longest runtime, but there was no significant difference in the runtime among all methods.




## Introduction

The Riemann proble [1,2] in basic gasdynamics refers to a setting where two adjacent regions of a compressible gas are initially separated by some barriers. Each region has a uniform set of conditions: density, velocity, and pressure. Suddenly, the barriers are removed and the energy imbalance between the two sides causes three waves such that some equilibrium average conditions are encouraged. These waves can be a shock wave (normal shock), a contact discontinuity (contact surface), or an expansion fan (rarefaction fan). It is possible to have two shocks with no fans, two fans with no shocks, and one shock and one fan. In either case, a contact discontinuity would be separating the two other waves. There is also a possibility of reaching a vacuum state in a zone surrounded by two contacts (one to the left and one to the right), which in turn separate two fans running away from each other [3]. However, the presence of vacuum is excluded in our work.

While the Riemann problem can be extended to more complicated environments [4,5], the simple Riemann problem is governed solely by the Euler equations (a system of time-dependent nonlinear hyperbolic partial differential equations,


*E-mail addresses:* osama.m@uob.edu.om, omarzouk@vt.edu







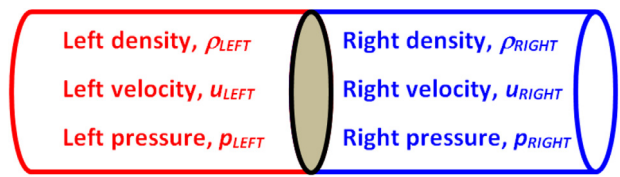

**Fig. 1.** Illustration of the Riemann problem. For a shock tube, the initial velocities are zeros.

**Table 1**
Initial left and right states in the Sod problem, with SI units.

| | Primitive variable | | | Derived quantities (with $\gamma = 1.4$) | | |
|---|---|---|---|---|---|---|
| | Density ($\rho$), kg/m³ | Velocity ($u$), m/s | Pressure ($p$), N/m² | Speed of sound ($a$), m/s | Specific internal energy ($e$), m²/s² | Specific enthalpy ($h$), m²/s² |
| Left state | 1.0 | 0.0 | 1.0 | 1.18322 | 2.5 | 3.5 |
| Right state | 0.125 | 0.0 | 0.1 | 1.05830 | 2.0 | 2.8 |

see Appendix A), which admit discontinuities in the solution. For a Computational Fluid Dynamics (CFD) solver implementing the Finite Volume Method (FVM), that is designed to handle compressible flows, the Riemann problem can be a good choice for testing the capability of the solver in forming fluxes at the finite volume face. The problem is focused on the convective/pressure flux in the one-dimensional space, thus eliminating complications of multi-dimensionality, inaccuracy of representing non-uniform boundary or initial conditions, and interference from other phenomena such as viscous effects and boundary layers, chemical reaction, mixing, multi-phase flow, heat transfer, thermal radiation, and source terms.

**Objectives and contribution**

Face flux construction is a core task in a finite volume solver for fluid flows as a simulation package. An end user of such a solver may be offered a large number of flux methods to select from for numerically handling convection (divergence) terms [6]. A good background in flux construction methods as well as an elementary analysis of their performance for a simple problem having a known analytical solution may help in making a good selection of a flux method in relevant simulations. To this end, the present work solves a particular example of the basic Riemann problem, known as the Sod problem with different numerical methods for building a flux vector at cell interfaces, and provides assessment of these flux methods.

The Riemann problem considered in this study is viewed as a good choice where the physical-chemical setting of the problem is simple due to the absence of the following features: variable gas properties, turbulence modeling, fluid-solid interaction, thermal radiation, chemical reactions, and multiphase flow. Also, the duration of the simulation time ensures that none of the traveling waves reaches the domain boundaries. Thus, the issue of wave reflection is not encountered. All these simplifications help emphasizing the role of the face flux method without much fear of being blurred due to treatments of several physical phenomena in the problem. On the other hand, the Riemann problem studied here contains both gradual changes and discontinuities, thus providing a level of challenge that can distinguish a good flux method. An overly simple test case with smooth variations everywhere does not offer a sufficient challenge for testing performance.

The contributions of the present work can be summarized as follows:

- Providing a short review of 22 face flux construction methods for the Euler equations, with a brief overview of each one which is sometimes accompanied by simple mathematical description
- Providing an assessment of these face flux construction methods in terms of the ability to handle a discontinuity, overall accuracy, and computational demand. Thus, helping a modeler and a user of a simulation kit in making a proper selection
- Providing quantitative details of the Sod problem (exact numerical values of some flow variables in each zone), helping a reader who wants to use it as a test problem to compare the values obtained by another solver with those reported here

**Sod's problem**

A shock tube is a special case of the Riemann problem in which the compressible gas on either side is initially stationary [7]. Fig. 1 gives a graphical illustration of the Riemann problem at the initial state. The Sod problem [8] is a shock tube problem with certain initial left and right densities and pressures, causing a right-going shock wave, a left-going expansion fan, and a right-going contact discontinuity separating the shock wave and the expansion fan.

Table 1 contrasts the left and right conditions in the Sod problem at the initial state (**initial time**: $t = 0$). The initial distribution of the primitive variables is demonstrated in Fig. 2. The jump is located at the middle position $x = 0.5$, which is midway between the left end of the problem domain ($x = 0$) and its right end ($x = 1$). Table 2 gives a summary of the fan-contact-shock wave pattern in the exact solution of the Sod problem.





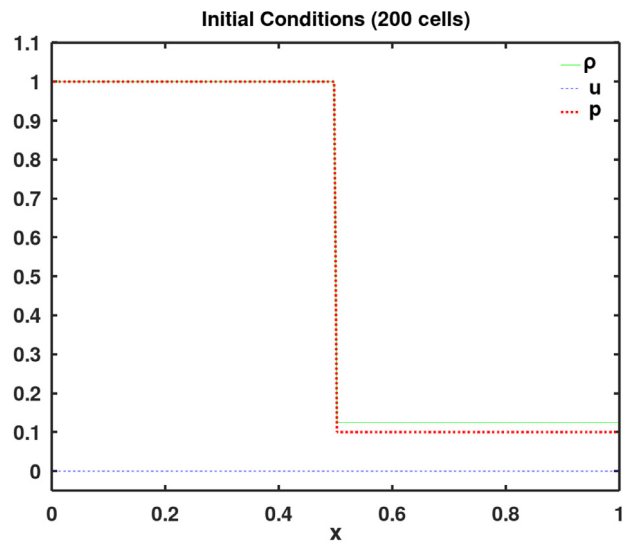

**Fig. 2.** Initial primitive variables in the Sod problem.

**Table 2**
Calculated properties of the three waves in the exact solution of the Sod problem.

| Wave type | Property | Value | Remarks |
|---|---|---|---|
| Expansion Fan | Velocity of the head | −1.18322 | The negative signs mean a leftward direction |
| | Velocity of the tail | −0.07027 | |
| Contact Discontinuity | Velocity | 0.92745 | The gas on both sides has the same velocity |
| | Pressure | 0.30313 | Equal pressures on both sides of the contact wave |
| | Left density | 0.42632 | $a = 0.99773$ |
| | | | $e = 1.77760$ |
| | | | $h = 2.48864$ |
| | Right density | 0.26557 | $a = 1.26411$ |
| | | | $e = 2.85354$ |
| | | | $h = 3.99496$ |
| Shock Wave | Velocity | 1.75216 | |
| | Shock-relative Mach number of unshocked gas | 1.65563 | Unshocked gas must be in a supersonic regime (Mach magnitude > 1) and the shocked has must be in a subsonic regime (Mach magnitude < 1) in order to |
| | Shock-relative Mach number of shocked gas | 0.65240 | avoid a universal entropy decrease, which violates the second law of thermodynamics [10] |

As a remark, the problem solved here reflects a virtual set-up. The values of the gas properties given here represent the physical ones including proper units. Using the international system (SI) of units [9], a density value of 1 means 1 kg/m³, and a pressure value of 1 means 1 Pa (pascal) or 1 N/m² (newton per square meter), and a velocity value of 0 means 0 m/s. Consequently, the speed of sound is in m/s, and the specific internal energy and the specific enthalpy are in m²/s² or J/kg (joules per kilogram). The equations solved are dimensional (no normalization or special scaling is done). With a goal of testing computational methods, the issue of whether or not this set-up can be realized experimentally may be of little concern.

At a simulation time of $t = 0.2$, the exact solution obtained for the Sod problem here is illustrated in Fig. 3, expressed in terms of the primitive variables and the specific internal energy. The abrupt change in the pressure designates the location of the shock wave traveling to the right. The abrupt change in the density that not accompanied by a change in pressure or velocity designates the location of the contact discontinuity, which also travels to the right. The gradual change in all the variables designates the extent of the expansion fan, which travels to the left. The left front of the fan is called "head", and the right edge of the fan is called "tail".

The computation of the exact solution presented earlier in this work, as well as the numerical ones to be presented later, was performed using computer codes written in GNU Octave, a scientific programming language [11]. The one-dimensional grid was in the form of 200 uniformly spaced cells.

As our exact solution is to be used as a reference for judging the numerical solution, it is important to check its validity. This is achieved here in two ways.

The first validity check of our implementation of the exact Riemann solver is achieved through satisfying the Rankine–Hugoniot shock condition ([12–15]) for the mass. The condition is used to compute a shock speed ($S$) other than the one





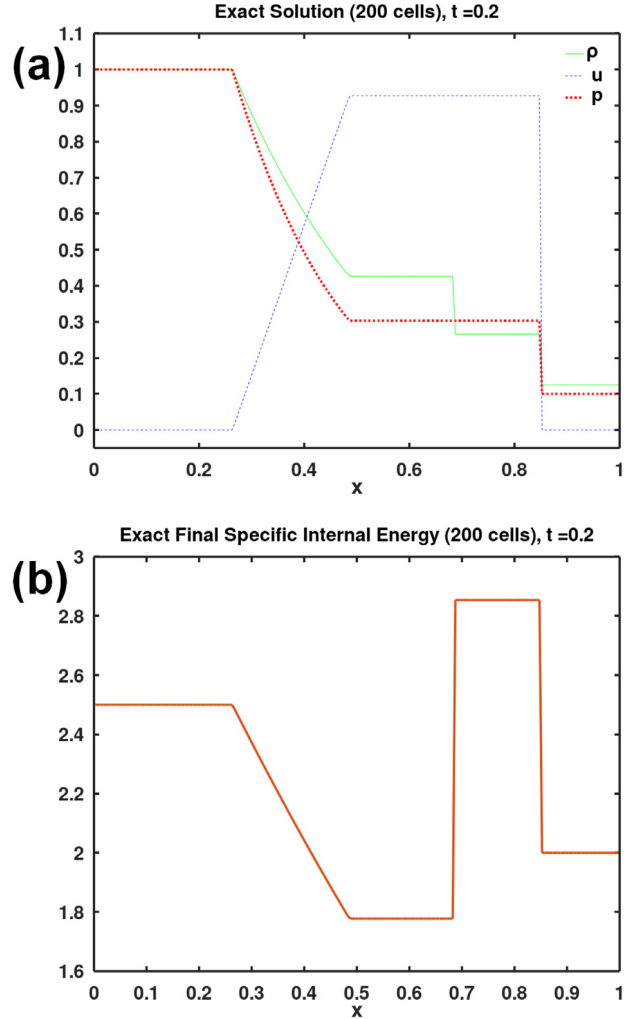

**Fig. 3.** Exact (analytical) solution obtained here at the final simulation time. Plot (a): primitive variables, density ($\rho$), velocity ($u$), and pressure ($p$). Plot (b): specific internal energy.

obtained from the exact Riemann solver (which was given already in Table 2). Based on the Rankine-Hugoniot mass condition, the shock wave speed should be

$$S = \frac{\Delta \rho u}{\Delta \rho} = \frac{(0.26557)\,(0.92745) - (0.125)(0)}{0.26557 - 0.125} = 1.75216 \tag{1}$$

where the operator $\Delta$ represents the difference across the shock wave. This calculated speed of 1.75216 matches the one given earlier in Table 2, which was found independently as

$$S = u_{RIGHT} + a_{RIGHT}\sqrt{1 + \frac{\gamma + 1}{2\,\gamma}\left(\frac{p_{CD}}{p_{RIGHT}} - 1\right)} \tag{2}$$

where $u_{RIGHT}$, $a_{RIGHT}$, and $p_{RIGHT}$ are the gas velocity, speed of sound, and pressure in the right undisturbed state; and $p_{CD}$ is the pressure around the contact discontinuity.

The second validity check of our implementation of the exact Riemann solver is achieved through a comparison with another independent solution given in Fig. 4. This independent solution is taken from the Castro project validation data [16], and Fig. 4 shows their results for the three primitive variables (density, velocity, and pressure) and the specific internal energy. Visual inspection suggests that the distributions of all these four quantities match those from solution given here in Fig. 3.





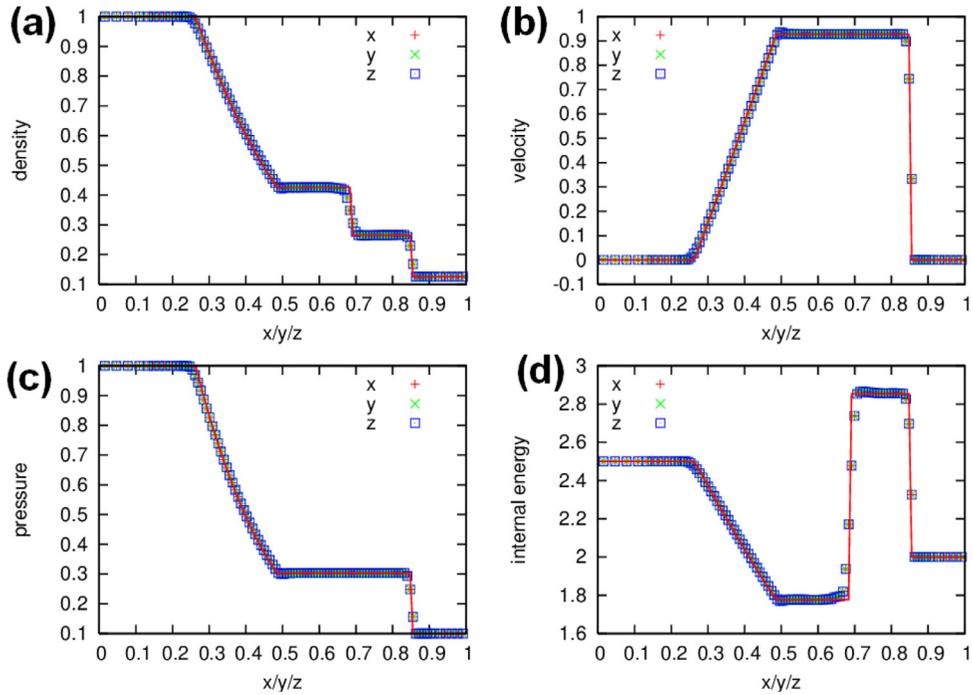

**Fig. 4.** Published exact and calculated final quantities: (a) density, (b) velocity, (c) pressure, and (d) specific internal energy using the Castor solver of compressible hydrodynamic equations for astrophysical flows at the same simulation time of ($t = 0.2$) as used in the present work. The exact solution is in a red solid line, whereas the calculated one is in markers (with courteous permission from Dr. Michael Zingale, Associate Professor, Stony Brook University, a member of the Castro team of core developers). The problem domain is three-dimensional. Source: https://amrex-astro.github.io/Castro/docs/Verification.html#sods-problem-and-other-shock-tube-problems.

## Face flux construction

The time integration of the conserved variables vector ($q$) is done using the Godunov scheme [17], which is a conservative explicit finite volume formula. For a uniformly spaced one-dimensional array of computational cells having a cell width of ($\Delta x$), let the cell center positon index be ($i$), then the left face flux vector is ($F_{i-0.5}$) and the right face flux vector is ($F_{i+0.5}$). Time integration of the conserved variables vector ($q$) over one time step ($\Delta t$) is performed as

$$q_i^{n+1} = q_i^n + \frac{\Delta t}{\Delta x} \left( F_{i+0.5}^n - F_{i-0.5}^n \right) \qquad (3)$$

The time index ($n$) refers to the latest known variable value, while the index ($n+1$) refers to the new value after time integration over one time step ($\Delta t$). The scheme is explicit, where the updating of a cell value is independent of updating other cells. This eliminates the need for solving a system of algebraic equations. Godunov scheme is of first order accuracy in time. The spatial accuracy depends on how the face flux vectors are constructed. It can be first order (lowest resolution) or second order (higher resolution).

Godunov scheme appears simple and straightforward, provided that the face flux vectors have been computed and are made ready to the scheme to use.

The process of preparing face flux vectors at intercell faces is split into two steps:

- Step 1 (MUSCL): use neighbor cell-center values of the primitive variables ($\rho$, $u$, $p$) to extrapolate and/or interpolate a 'fictitious' value at both sides of the facial interface. The left side of the face has the fictitious primitive values denoted by: $\rho_L$, $u_L$, $p_L$. These are the local left values. Similarly, the right fictitious primitive values of the face are denoted by: $\rho_R$, $u_R$, $p_R$. These are the local right values. It should be noted that although physically these locally left and locally right values are defined at the same spatial location, they are not necessarily equal. Once these fictitious face-left and face-right values are computed, a local fictitious Riemann problem is constructed. There is a jump in one or more local flow variables across the face. Along with the interpolation/extrapolation, the preparation of the local Riemann problem at each face includes also the use of a flux limiter (also called a slope limiter) through a scheme called MUSCL, which stands for Monotonic Upstream-centered Scheme for Conservation Laws [18]. This scheme is a nonlinearized version of the so-called $\kappa$ 3-point scheme of parabolic interpolation/extrapolation [19]. Whereas the interpolation/extrapolation coefficients in the $\kappa$ scheme are solution-independent, the MUSCL scheme adds to the $\kappa$ scheme a solution-dependent limiter, which downgrades the interpolation/extrapolation level near sharp changes (like a shock wave) to use only the





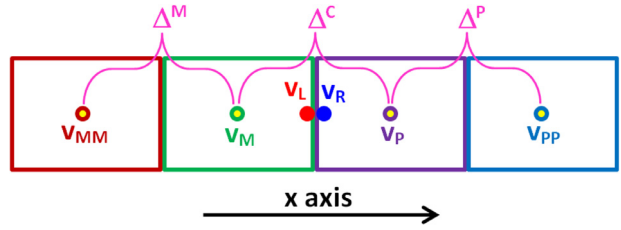

**Fig. 5.** Illustration of interpolation/extrapolation process to obtain the left and right face values of a generic variable ($v$) using the MUSCL scheme.

nearest cell-center value, but upgrades the interpolation/extrapolation level to allow using more cell-center values in regions of smooth variations. This adaptive behavior allows the use of the high-resolution with its advantage of spatial accuracy without suffering from severe oscillations near sharp changes by 'trimming' or even 'switching off' the high-resolution capability and keeping only the first-order spatial accuracy. There are different flux limiters [21,22]. Here, the van Leer flux limiter [20] is used, which has intermediate characteristics among others in terms of the dissipation level [23].

- Step 2 (flux method): After a local fictitious Riemann problem is constructed, then comes the step of inferring a face flux vector based on the jump across the local left and local right flow variables. This step is the core of this work. The work considers 22 methods of face flux methods proposed in the literature, and compare their performance for the Sod problem as a test case.

To explain the MUSCL step, consider a generic variable ($v$) that has four values defined at cell centers of four consecutive cells labeled in the direction of the increasing x-coordinate as *MM*, *M*, *P*, and *PP*. The cell values (cell center values) are labeled in the direction of the increasing x-coordinate as $v_{MM}$, $v_M$, $v_P$, and $v_{PP}$ according to in Fig. 5. The goal is to construct a limited interpolated/extrapolated fictitious left and right values at the face separating the middle cells *M* and *P*.

Let

$$\Delta^M = v_M - v_{MM} \tag{4}$$

$$\Delta^C = v_P - v_M \tag{5}$$

$$\Delta^P = v_{PP} - v_P \tag{6}$$

Two consecutive gradient ratios are calculated, left one and right one, respectively, as

$$r_L = \frac{\Delta^C}{\Delta^M} \tag{7}$$

$$r_R = \frac{\Delta^C}{\Delta^P} \tag{8}$$

There is a possibility of a division by zero when calculating $r_L$ and $r_R$. However, this is avoided by applying the following condition

$$\text{if } \left|\Delta^M\right| \leq \epsilon, \quad r_L = 0 \tag{9}$$

$$\text{if } \left|\Delta^P\right| \leq \epsilon, \quad r_R = 0 \tag{10}$$

where $\epsilon$ is a small threshold, which was decided by the computing program as $2.22 \times 10^{-16}$.

Then, two flux limiting values $\phi_L$ and $\phi_R$ are obtained for the left and right sides of the face, respectively, using the van Leer flux limiter function.

$$\phi_L = \frac{r_L + |r_L|}{1 + |r_L|} \tag{11}$$

$$\phi_R = \frac{r_R + |r_R|}{1 + |r_R|} \tag{12}$$





Now, the final stage is to compute the fictitious left and right face values $v_L$ and $v_R$, respectively, using the following version of the MUSCL scheme:

$$v_L = v_M + 0.5 \; \phi_L \; \Delta^M \tag{13}$$

$$v_R = v_P - 0.5 \; \phi_R \; \Delta^P \tag{14}$$

In the case of unity flux limiters ($\phi_L = \phi_R = 1$), the above equations represent a pseudo Linear-Upwind Interpolation (LUI), which is also called Second Order Upwind (SOU). In that case, each of the fictitious face values is constructed by a standard linear extrapolation from the two nearest cell-center values on the same face side (either left or right). This pure extrapolation (no interpolation) feature corresponds to setting the parameter value $\kappa = -1$ in the linear version of the MUSCL scheme (the $\kappa$ scheme).

On the other hand, if the flux limiters are set to zero ($\phi_L = \phi_R = 0$), then a pseudo first order (classical upwinding) interpolation is obtained, where each fictitious face value is simply taken as the nearest cell-center value on the same face side (either left or right). The term 'pseudo' is used in the above text because in true upwinding schemes, taking the left or right cell-center value depends on the speed of the characteristic [24], also called the wave speed or the information speed, relevant to the concerned variable at the face. Such checking of the direction of the information is not done here, and it is always assumed that for the left fictitious face value, the flow information (the wind) comes from the left, and for the right fictitious face value, it comes from the right.

### Analyzed flux methods

Twenty two flux methods are tested in this work and they are listed here with a brief description on each one. Giving detailed explanation or the mathematical formulation behind each method can make this section extremely long. Instead, the interested reader is cordially requested to visit the literature for more information about any of methods found of interest. In presenting the data, a short name (a code) is assigned to each method as a quick unique identifier. The methods are listed below by their codes. The sequence of introducing the methods here is the same one used later when presenting their graphical results and the detailed root mean square errors of their numerical results.

*Riemann*

This flux method is also called Godunov solver. Godunov is credited to formulate this exact (analytic) Riemann solver in the same reference cited before when addressing the Godunov time integration scheme. This flux method solves a Riemann problem exactly at each intercell face, with the local fictitious face-left and face-right conditions provided by the MUSCL scheme replacing the global initial physical left and right states of the gas in the Sod problem separated by the middle barriers. This local Riemann solver identifies the state of the gas at the face interface, such as a zone between a left-going expansion fan and a right-going contact discontinuity, a zone between a right-going shock wave and a right-going contact discontinuity, or a zone within an expansion fan whose tail and head are going opposite to each other, which is called a transonic fan as it accelerates the gas from subsonic (Mach number magnitude < 1) to supersonic (Mach number magnitude > 1). Then, the face flux is constructed using the primitive values of the identified state. Thus, if the identified state has numerical values of the density, velocity, and pressure denoted by $\rho_0, u_0$, and $p_0$, then the constructed face flux is simply

$$F_{Face} = \left\{ \begin{array}{c} \rho_0 \; u_0 \\ \rho_0 u_0^2 + p_0 \\ \left( \frac{\gamma p_0}{(\gamma - 1)} + \frac{1}{2} \rho_0 u_0^2 \right) u_0 \end{array} \right\}$$

*Roe*

This is the Roe approximate solver for face flux, which is also known as Roe upwind scheme [25]. The nonlinear flux vector in Euler equations is replaced by a locally linearized approximate flux function, which utilizes a $3 \times 3$ matrix called a Roe-average Jacobian matrix. The Roe flux method gives three equally-spaced waves. It cannot capture the spread of the expansion fan. Instead, it represents it as a zero-thickness wave. Denoting the MUSCL returned face-left values of the density and velocity by $\rho_L$ and $u_L$, respectively; and the MUSCL returned face-right values of the density and velocity by $\rho_R$ and $u_R$, respectively; then, the Roe method introduces a customized square-root-density weighted average velocity defined as: $u_{Roe} = (u_L \sqrt{\rho_L} + u_R \sqrt{\rho_R})/(\sqrt{\rho_L} + \sqrt{\rho_R})$. Let the face-left specific total enthalpy be $h_{T,L} = \gamma \, p_L/(\rho_L(\gamma - 1)) + 0.5 \, u_L^2$, where $p_L$ is the face-left pressure as provided by the MUSCL scheme. Similarly, let the face-right specific total enthalpy be $h_{T,R} = \gamma \, p_R/(\rho_R(\gamma - 1)) + 0.5 \, u_R^2$, where $p_R$ is the face-right pressure as provided by the MUSCL scheme. The Roe method introduces a square-root-density weighted average total specific enthalpy defined as: $h_{T,Roe} = (h_{T,L} \sqrt{\rho_L} + h_{T,R} \sqrt{\rho_R})/(\sqrt{\rho_L} + \sqrt{\rho_R})$. Finally, a Roe-specific speed of sound is defined as: $a_{Roe} = \sqrt{(\gamma - 1)(h_{T,Roe} - 0.5 \, u_{Roe}^2)}$. The three wave speeds in Roe method





(also called Roe eigenvalues) are: $u_{Roe}$, $u_{Roe} - a_{Roe}$, and $u_{Roe} + a_{Roe}$. The three Roe specific quantities ($u_{Roe}$, $h_{T,Roe}$, and $a_{Roe}$) are used to construct the face flux vector.

*KNP*

This is the central-upwind method of Kurganov, Noelle, and Petrova for the face flux [26]. The term central-upwind here comes from the notice that the method is central in principle (thereby, offering high resolution with small dissipation), but still has an upwinding character by respecting the direction of the wave travel. The method does not attempt to approximate the exact flux as in the case of the Roe method.

*KT*

This is another central-upwind method, developed by Kurganov and Tadmor for forming the intercell face flux vector [27]. It resembles the KNP method in the structure, with the KT method (which is chronologically earlier) being simpler. For example, the KT method assigns fixed weights of 0.5 to the face-left and face-right MUSCL flux vectors, whereas the KNP method uses flow-based weights.

*SW*

This is the Steger and Warming Flux Vector Splitting (FVS) method for face flux [28,29]. The face flux vector is computed as a sum of two flux vectors: a 'plus' flux vector and a 'minus' flux vector, as $F = F^+ + F^-$. The subvector $F^+$ is associated with the positive eigenvalues of the Jacobian matrix $\partial F / \partial q$, while the subvector $F^-$ is associated with its negative eigenvalue.

*vanLeer*

This is another Flux Vector Splitting (FVS) method, developed by van Leer [30]. It allows a continuous slope of two formulation functions at the special values of Mach numbers of $-1$, 0, and 1. A quadratic function of the fictitious face-left and face-right Mach numbers is introduced. This addresses a shortcoming in the chronologically earlier SW Flux Vector Splitting method.

*AUSM*

This is the Advection Upstream Splitting Method (AUSM) for face flux construction [31]. The flux is divided into a pressure component and a convective component for better treatment.

*AUSM$^+$*

This is a follow-up improvement [32] of the Advection Upstream Splitting Method (AUSM) for computing face fluxes by addressing observed shortcomings in some cases in which the original AUSM technique did not maintain its tradition of accuracy and efficiency.

*AUSM$^{+-}$up (for all speeds)*

This is an extension of AUSM$^+$ to have the performance nearly independent of the Mach number for the low-Mach subsonic limit up to a freestream Mach number of 0.5. This is achieved through proper scaling using asymptotic analysis [33]. Also the performance level is nearly comparable in high-Mach subsonic flows and supersonic flows.

*AUFS*

AUFS stands for Artificially Upstream Flux Vector Splitting method [34]. It employs two simultaneous Flux Vector Splitting mechanisms, one by dividing the flux vector into pressure and non-pressure components, and another by using a method-specific Mach number (one that is based on the simple average of the face-left and face-right velocities, $u_{avg} = (u_L + u_R)/2$) as a splitting weight. The method also respects the sign of the average velocity of the face-left and face-right fictitious values. The method is relatively easy to implement.

*HLL-Davis1*

This is the Harten, Lax, van Leer (HLL) 2-wave model [35], with Davis first estimation of right and left wave speeds [36]. Let $\rho_L$, $u_L$, and $p_L$ be the face-left fictitious values of the density, velocity, and pressure, respectively; and $\rho_R$, $u_R$, and $p_R$ be their face-right counterparts. Then, a face-left speed of sound is calculated as $a_L = \sqrt{\gamma \ p_L / \rho_L}$, and face-right speed of sound is calculated as $a_R = \sqrt{\gamma \ p_R / \rho_R}$. This version of the HLL family of methods sets the left wave speed as $S_L = u_L - a_L$, and the right wave speed as $S_R = u_R + a_R$.





*HLL-Davis2*

This is a second variant of the Harten, Lax, van Leer (HLL) 2-wave model, with Davis second estimation of right and left wave speeds. With the same meaning of $u_L, a_L, u_R, a_R$ as used in the HLL-Davis1 method earlier, this HLL-Davis2 version of the HLL family sets the left wave speed as $S_L = \min(u_L - a_L, u_R - a_R)$, and the right wave speed as $S_R = \max(u_L + a_L, u_R + a_R)$. The min( ) function gives the smaller of its two inputs, and the max( ) function gives the larger of them.

*HLL-Roe*

This is a third variant of the Harten, Lax, van Leer (HLL) 2-wave model, with two of the Roe eigenvalues used for the right and left wave speeds. Using same notations in the Roe method and in the HLL-Davis1 method described earlier, the two wave speeds in the HLL-Roe method described now are: $S_L = u_{Roe} - a_{Roe}$ and $S_R = u_{Roe} + a_{Roe}$.

*HLL-Einfelt*

This is a fourth variant of the Harten, Lax, van Leer (HLL) 2-wave model, with Einfelt eigenvalues for the right and left wave speeds [37]. It is also called (HLLE) flux method. The left and right wave speeds here are: $S_L = u_{Roe} - u_{Einfelt}$ and $S_R = u_{Roe} + u_{Einfelt}$, respectively. The symbol $u_{Roe}$ refers to the Roe square-root-density weighted averaged velocity as described in the Roe method. The velocity $u_{Einfelt}$ is computed as: $u_{Einfelt} = \sqrt{(a_L^2\sqrt{\rho_L} + a_R^2\sqrt{\rho_R})/(\sqrt{\rho_L} + \sqrt{\rho_R}) + 0.5\sqrt{\rho_L}\sqrt{\rho_R}(u_R - u_L)^2/(\sqrt{\rho_L} + \sqrt{\rho_R})^2}$. A subscript $L$ refers to a face-left value, and a subscript $R$ refers to a face-right value. The face-left speed of sound $a_L$ and the face-right speed of sound $a_R$ are as discussed in the HLL-Davis1 method.

*HLL-pBased*

This is a fifth variant of the Harten, Lax, van Leer (HLL) 2-wave model, with pressure-based calculation for the right and left wave speeds [38]. The method introduces a customized threshold or criterion pressure $p^*$ used to determine the left and right wave speeds, where $p^* = 0.5(p_L + p_R) + 0.125(u_R - u_L)(a_R + a_L)(\rho_R + \rho_L)$. The left and right wave speeds here are: $S_L = u_L - f_L a_L$ and $S_R = u_R + f_R a_R$, respectively. The symbols in these formulas have the same meanings as those given for the HLL-Davis1 method. Two new symbols appear here ($f_L$ and $f_R$), which are factors decided based on $p^*$ and how it compares with the face-left pressure $p_L$ and the face-right pressure $p_R$, as follows: if $p^* \le p_L$, $f_L = 1$, if $p^* > p_L$, $f_L = \sqrt{1 + (p^*/p_L - 1)(\gamma + 1)/(2\gamma)}$; if $p^* \le p_R$, $f_R = 1$, and if $p^* > p_R$, $f_R = \sqrt{1 + (p^*/p_R - 1)(\gamma + 1)/(2\gamma)}$.

*HLLC-Davis1*

This is an extended 3-wave version [39] of the Harten, Lax, van Leer, Davis (HLL-Davis1) 2-wave model, where a third wave (contact discontinuity or contact surface) is included. The additional "C" letter in the method code stands for the added "contact" wave. This method for flux construction is comparable in accuracy and robustness to that of the exact Riemann flux method, while being simpler and computationally less demanding than it. Davis first estimation of right and left wave speeds is used.

*HLLC-Davis2*

This is an extended 3-wave version of the Harten, Lax, van Leer, Davis (HLL-Davis2) 2-wave model, where a third wave (contact discontinuity) is added. Davis second estimation of right and left wave speeds is used.

*HLLC-Roe*

This is an extended 3-wave version of the Harten, Lax, van Leer, Roe (HLL-Roe) 2-wave flux model, where a third wave (contact discontinuity) is added.

*HLLC-Einfelt*

This is an extended 3-wave version of the Harten, Lax, van Leer, Einfelt (HLL-Einfelt or HLLE) 2-wave flux model, where a third wave (contact discontinuity) is added.

*HLLC-pBased*

This is an extended 3-wave version of the Harten, Lax, van Leer, pressure-based (HLL-pBased) 2-wave flux model, where a third wave (contact discontinuity) is added.





*LF*

The Lax-Friedrichs (LF) method, which is also called the Lax method [40] is a basic explicit three-point technique used for solving time-dependent partial differential equations. It is first-order accurate in space. This method gives the most dissipative flux among all basic stable explicit flux construction methods. The method does not require solving a Riemann problem (either exact or approximate) at each cell face. The face flux is constructed as $F_{face} = 0.5 \ (F_L + F_R) - 0.5\Delta x/\Delta t(q_R - q_L)$, where $F_L$ and $F_R$ are the flux vectors after numerical evaluation using the face-left and face-right values of the primitive variables, respectively; $q_L$ and $q_R$ are the conserved variables after numerical evaluation using the face-left and face-right values of the primitive variables, respectively; $\Delta x$ is the cell width (the spatial step); and $\Delta t$ is the time step (the temporal increment).

*Rusanov*

This is a single-wave simplification of the Roe 3-wave model [41]. It is also called generalized Lax-Friedrichs flux method. Like the Lax-Friedrichs flux method, this method does not require solving a Riemann problem (either exact or approximate) at each cell face. However, unlike the Lax-Friedrichs flux method, it makes use of the local speed of sound. The flux calculation formula is nearly identical to the one given for the Lax-Friedrichs flux method, with the only difference is that the uniform virtual speed ($\Delta x/\Delta t$) is replaced by a gas-dependent speed $\lambda_{max} = \max(|u_L| + a_L, \ |u_R| + a_R)$, where $|u_L|$ and $|u_R|$ are the absolute values of the face-left and face-right gas velocities, respectively; and $a_L$ and $a_R$ are the speeds of sounds obtained using the face-left and face-right densities and pressures, respectively, as mentioned for the HLL-Davis1 method.

## Results

The Sod problem is solved numerically using each of the listed 22 flux methods in the previous section. No other change was done from one numerical simulation to another. Common simulation settings include a cell size of $\Delta x = 0.005$ (thus, the shock tube of unit length is divided into 200 equal cells), and a time step of $\Delta t = 0.001$, which was related to the cell size as

$$\Delta t = \frac{\text{Co}_{max} \ \Delta x}{S_{max}} \tag{1}$$

where $\text{Co}_{max}$ is the desired upper limit of the Courant number (Courant–Friedrichs–Lewy number, CFL number) as a condition for spatio-temporal stability in explicit time integration of partial differential equations, requiring the Courant number (Co) to be below unity. As the Co decreases, stabilizing numerical dissipation increases [42]. The parameter $S_{max}$ is an estimate of the maximum wave speed in the problem. In the present work, $S_{max} = 2$, which is a reasonable value for the highest possible sum of the local absolute gas velocity and the speed of sound in the problem, $|u| + a$. The target maximum Courant number was set to $\text{Co}_{max} = 0.4$. This target was largely satisfied in all performed 22 simulations, as found by post-processing the simulation results.

In Figs. 6, 7, 8, 9, 10, 12, 13, 14 and 15, the spatial distribution of the three primitive variables (density $\rho$, gas velocity $u$, and pressure $p$) at the final simulation time of 0.2 (thus, after 200 time steps) are given. The title of each figure includes

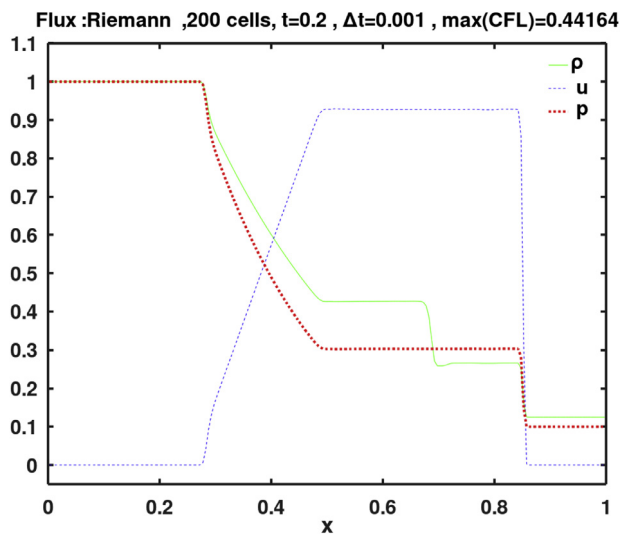

**Fig. 6.** Calculated final primitive variables using face flux of exact Riemann solver (also called Godunov flux).





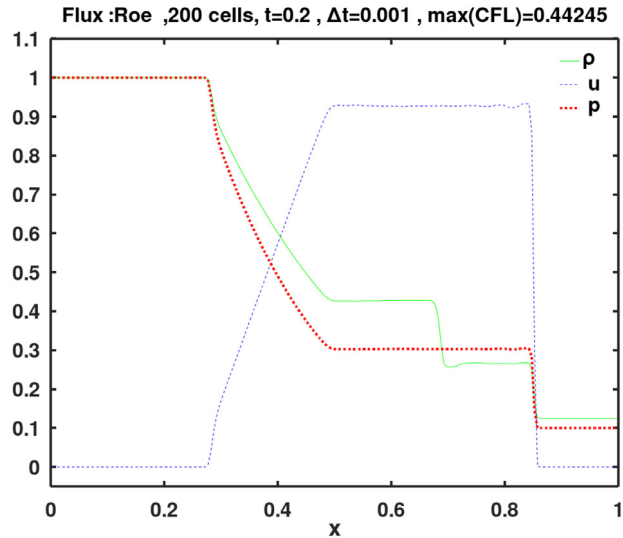

**Fig. 7.** Calculated final primitive variables using face flux of Roe.

the code (the short name) of the flux method used in the simulation. In addition, the maximum observed Co number (CFL number) is also listed in the plots titles. The order of presenting these figures for the different flux methods is the same order used when briefly describing them earlier.

The LF method in Fig. 14 shows significant undesirable dissipation causing severe smearing of jumps in the solution. No other method showed such level of unacceptable performance.

For the AUSM family (Fig. 10), the basic method exhibits a small overshoot at the tail of the fan. AUSM$^+$ did not resolve this problem, but it disappears in the further-improved version of AUSM$^{+-}$up.

The flux methods of Roe, KT, SW, AUFS, HLL-Davis1, HLL-Roe, HLL-Einfelt, HLL-pBased, HLLC-Davis1, HLLC-Roe, HLLC-Einfelt, HLLC-pBased, and Rusanov show small velocity oscillations (wiggles) to the left of the shock (i.e., after the gas is shocked). These wiggles are particularly noticeable with the SW method.

The flux methods of the exact Riemann solver (Fig. 6), KNP (Fig. 8a), vanLeer (Fig. 9b), AUSM$^{+-}$up (Fig. 10c), HLL-Davis2 (Fig. 12b), and HLLC-Davis2 (Fig. 13b) appear to perform the best in terms of being free from visual oscillations in the solution.

Within the expansion fan, near its head, the methods of AUSM and AUSM$^+$ (Fig. 10), and AUFS (Fig. 11) do not capture well the spatial increase of the gas velocity. They produce a fan head that is shifted to the right of its exact location. On the other hand, the methods of KT (Fig. 8b), and SW and vanLeer (Fig. 9) satisfactorily predict the propagation of the fan head.

To quantify the closeness of each numerical solution to the exact one, the root mean square error (RMSE) for each of the three primitive variables was computed. Let ($v$) denote any of these variables. At each cell-center, let the exact value be ($v_{ex}$) while the obtained one by numerically solving Euler equations using a certain flux method be ($v_{num}$). Then the RMSE for the generic variable ($v$) is

$$\text{RMSE(v)} = \sqrt{\frac{\sum_{i=1}^{n} \left(v_{num,i} - v_{ex,i}\right)^2}{n}} \tag{2}$$

where ($i$) is an index for the cell-centers, and ($n$) is the total number of cell-centers (which is 200 in this work).

The full list of these errors is given in Table 3. The flux methods are listed in the table in the same sequence of introducing them earlier. The velocity seems to be the most challenging variable to predict, because it always has the highest RMSE value among the three primitive variables regardless of the flux method.

The maximum magnitudes of the three promotive variables (density, gas velocity, and pressure) in the exact solution are of the same order of magnitude. In addition, their individual RMSE values are also of the same order of magnitude. It may be then helpful to lump their RMSE values into a single measure of inaccuracy by simple addition, leading to a single aggregate RMSE for each flux method. This is done in Table 3, by adding the total RMSE value for each flux method in the last column. If the total RMSE values are ordered ascendingly, the method of Roe comes first with the smallest total deviation from the exact solution (total RMSE = 0.03789). For all flux methods other than AUSM$^{+-}$up, AUSM, SW, and LF, the total RMSE increases gradually with small increments from a method to the next less-accurate one up to AUSM$^+$ (total





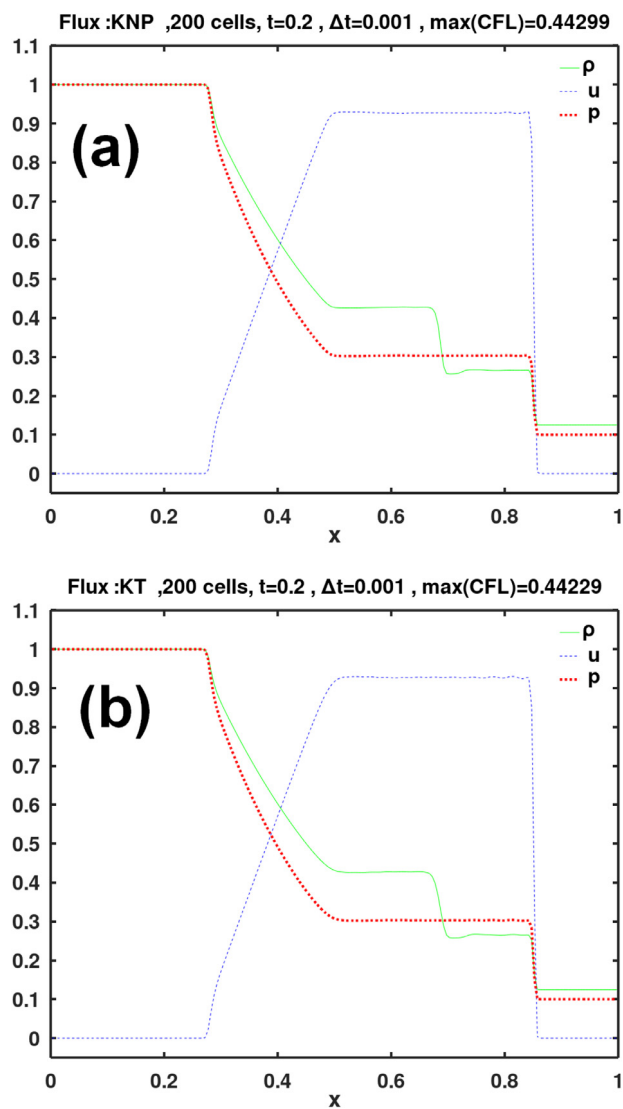

**Fig. 8.** Calculated final primitive variables using face flux of the central-upwind methods of (a) Kurganov, Noelle, and Petrova (KNP) and (b) Kurganov and Tadmor (KT).

RMSE = 0.04367). This is followed by mild inaccuracy increases for the method AUSM$^{+-}$up (total RMSE = 0.0449) and AUSM (total RMSE = 0.05037). Then, as the end of the ordered list, the total RMSE jumps by a factor of about 4 for the SW and LF methods, thereby clearly delineating their solutions as improper.

While the accuracy of a flux method is very important when assessing it, its computational demand is worthy of consideration. In large three-dimensional finite volume simulations, face flux construction is a process that happens millions of times to advance the solution by just one time step, and the solution may need thousands of time steps to complete. Therefore, even an extremely small time saving in the calculation of the flux vector at a face can have a useful impact when considering the total runtime of a realistic simulation. The elapsed clock time of the simulation excluding the processes of (1) initializing variables, (2) calculating an exact solution (as a reference to compare with), (3) plotting results, and (4) writing output data files was recorded for each flux method, and the results are presented in Table 4 in an ascending order. It is worthy of mentioning here that the relative time (one method relative to another) is more meaningful than the absolute time spent. The computing machine used in this work is a laptop Dell Inspiron (model N5010), with 6 GB of installed RAM memory and an Intel® Core™ i5 CPU, M 480 @ 2.67 GHz. The SW method took the shortest time of 41.13235s, while the Riemann exact solver took the longest time of 46.49466 s. The ratio of these two extremum values is 1.13037.





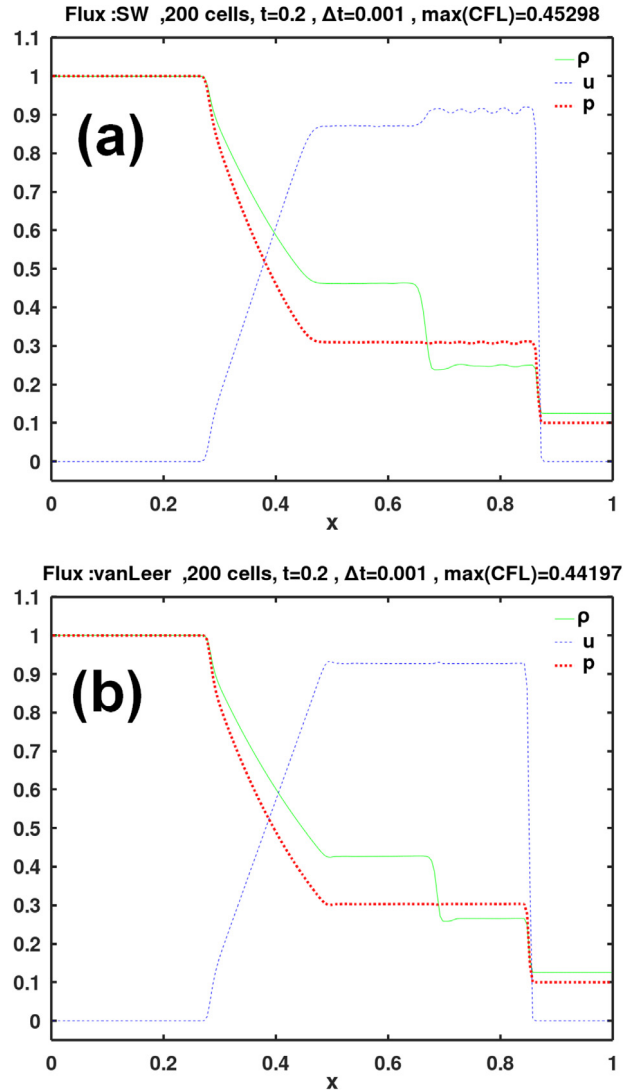

**Fig. 9.** Calculated final primitive variables using face flux of the Flux Vector Splitting methods of (a) Steger and Warming (SW) and (b) van Leer.

A percentage increase of the runtime (using the shortest recorded one as a reference) is included in Table 4. This relative measure of computational cost is calculated as

$$\% \text{ extra runtime} = \frac{\text{runtime (certain method)} - \text{ shortest run time}}{\text{shortest run time}} \times \ 100\% \tag{3}$$

The changes in the elapsed time are not large, with the extra time not exceeding 13% when using the exact Riemann solver flux method (as the one found to be the most computationally demanding) compared to any other method. This can be attributed to the inherent simplicity of the problem.

## Conclusions

This work provided a validated exact solution of the Sod problem, which is a one-dimensional jump problem for an inviscid flow, featuring three traveling waves in a tube: shock wave, contact discontinuity, and expansion fan. Quantitative details about the exact solution were provided, which can be useful as a reference to compare with. The exact solution also served as a tool to assess the performance of 22 flux methods used while generating an approximate computational solution for the same problem by solving the one-dimensional Euler equations. The difference in the accuracy and elapsed time among the obtained numerical solutions is related to the underlying flux method. The Godunov explicit finite volume scheme was used for time integration. A full extrapolation version of the MUSCL scheme was used along with the van Leer flux limiter.





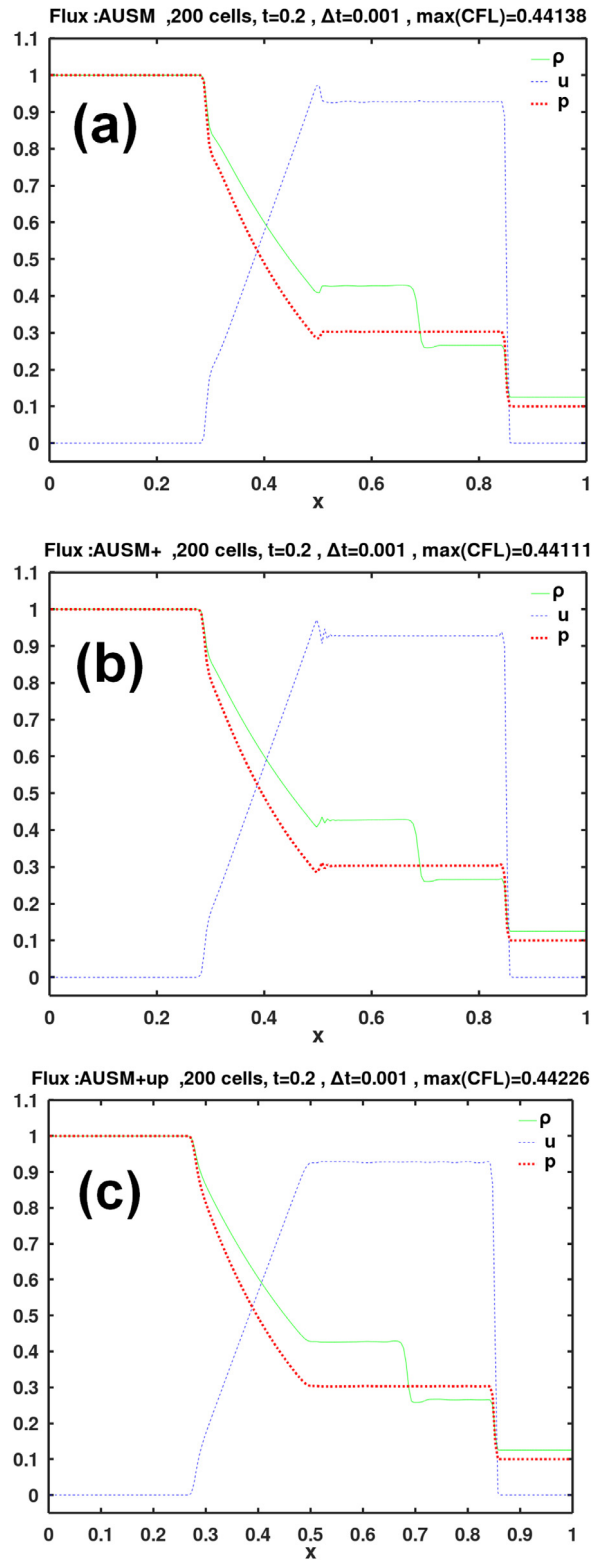

**Fig. 10.** Calculated final primitive variables using face flux of the Advection Upstream Splitting Method family of methods: (a) basic AUSM, (b) improved AUSM[+], and (c) second improved AUSM[+]−up, for all speeds.





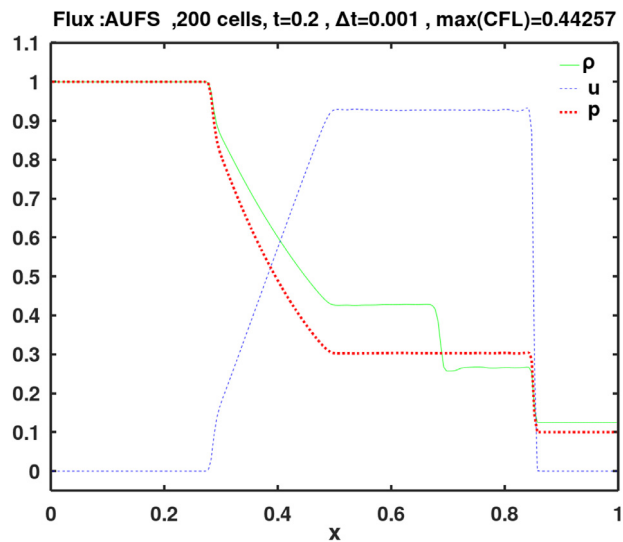

**Fig. 11.** Calculated final primitive variables using face flux of the Artificially Upstream Flux Vector Splitting (AUFS) method.

**Table 3**
Comparison of the RMSE (root mean square error) for each of the three primitive variables and their sum, corresponding to the flux methods used here.

| Index | Flux method (short name) | RMSE (density) | RMSE (velocity) | RMSE (pressure) | RMSE (total) |
|-------|--------------------------|----------------|-----------------|-----------------|--------------|
| 1 | Riemann | 0.00798 | 0.02345 | 0.00811 | 0.03954 |
| 2 | Roe | 0.00777 | 0.02216 | 0.00796 | 0.03789 |
| 3 | KNP | 0.00829 | 0.02423 | 0.00807 | 0.04059 |
| 4 | KT | 0.00889 | 0.02519 | 0.00760 | 0.04168 |
| 5 | SW | 0.03281 | 0.11764 | 0.02876 | 0.17921 |
| 6 | vanLeer | 0.00767 | 0.02624 | 0.00758 | 0.04149 |
| 7 | AUSM | 0.01127 | 0.02595 | 0.01315 | 0.05037 |
| 8 | AUSM+ | 0.00947 | 0.02380 | 0.01040 | 0.04367 |
| 9 | AUSM+-up | 0.00748 | 0.03047 | 0.00695 | 0.04490 |
| 10 | AUFS | 0.00835 | 0.02121 | 0.00891 | 0.03847 |
| 11 | HLL-Davis1 | 0.00818 | 0.02184 | 0.00788 | 0.03790 |
| 12 | HLL-Davis2 | 0.00829 | 0.02423 | 0.00807 | 0.04059 |
| 13 | HLL-Roe | 0.00821 | 0.02213 | 0.00796 | 0.03830 |
| 14 | HLL-Einfelt | 0.00821 | 0.02219 | 0.00797 | 0.03837 |
| 15 | HLL-pBased | 0.00824 | 0.02312 | 0.00799 | 0.03935 |
| 16 | HLLC-Davis1 | 0.00793 | 0.02234 | 0.00805 | 0.03832 |
| 17 | HLLC-Davis2 | 0.00790 | 0.02381 | 0.00800 | 0.03971 |
| 18 | HLLC-Roe | 0.00787 | 0.02209 | 0.00794 | 0.03790 |
| 19 | HLLC-Einfelt | 0.00788 | 0.02213 | 0.00794 | 0.03795 |
| 20 | HLLC-pBased | 0.00786 | 0.02324 | 0.00797 | 0.03907 |
| 21 | LF | 0.04383 | 0.11586 | 0.05071 | 0.21040 |
| 22 | Rusanov | 0.00889 | 0.02519 | 0.00760 | 0.04168 |

The analysis of results suggests that there is a number of flux construction methods that can offer good performance with a similar level of speed, including the methods of Kurganov-Noelle-Petrova, van Leer, Harten-Lax-van Leer-Davis2, and Harten-Lax-van Leer-Contact-Davis2. On the other hand, two flux methods may be avoided due to high oscillation or excessive dissipation.

In terms of accuracy (expressed as a small sum of root mean square errors for density, velocity, and pressure), the best flux method among all the 22 analyzed ones was found to be the Roe method. In terms of computational cost (represented by the runtime), the best flux method was the Steger and Warming Flux Vector Splitting method, and then comes the Advection Upstream Splitting Method (AUSM).





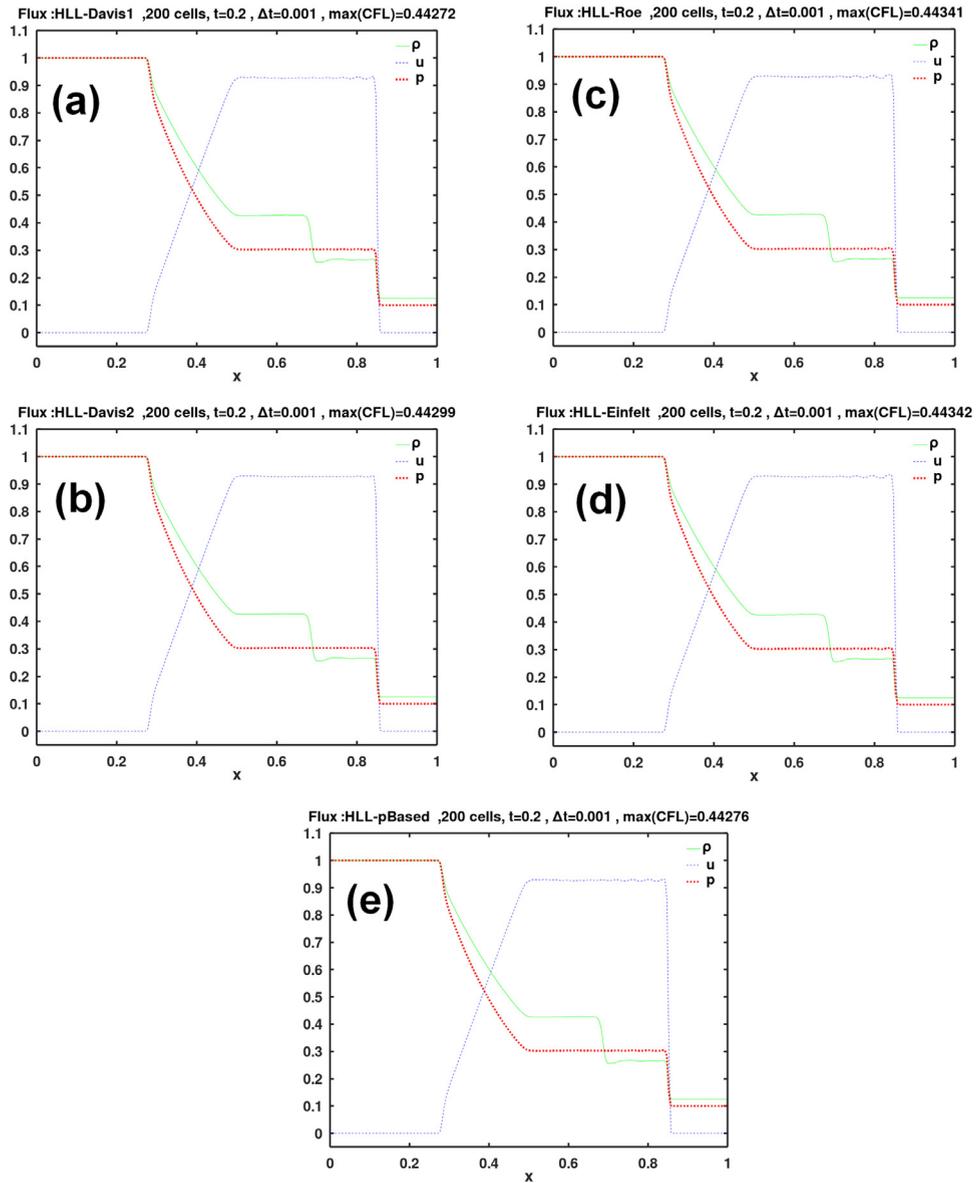

**Fig. 12.** Calculated final primitive variables using face flux of the Harten, Lax, van Leer (HLL) family of methods. The right/left wave speeds are obtained using: (a) first estimation of Davis, (b) second estimation of Davis, (c) eigenvalues of Reo, (d) eigenvalues of Einfelt, and (e) pressure-based calculation.





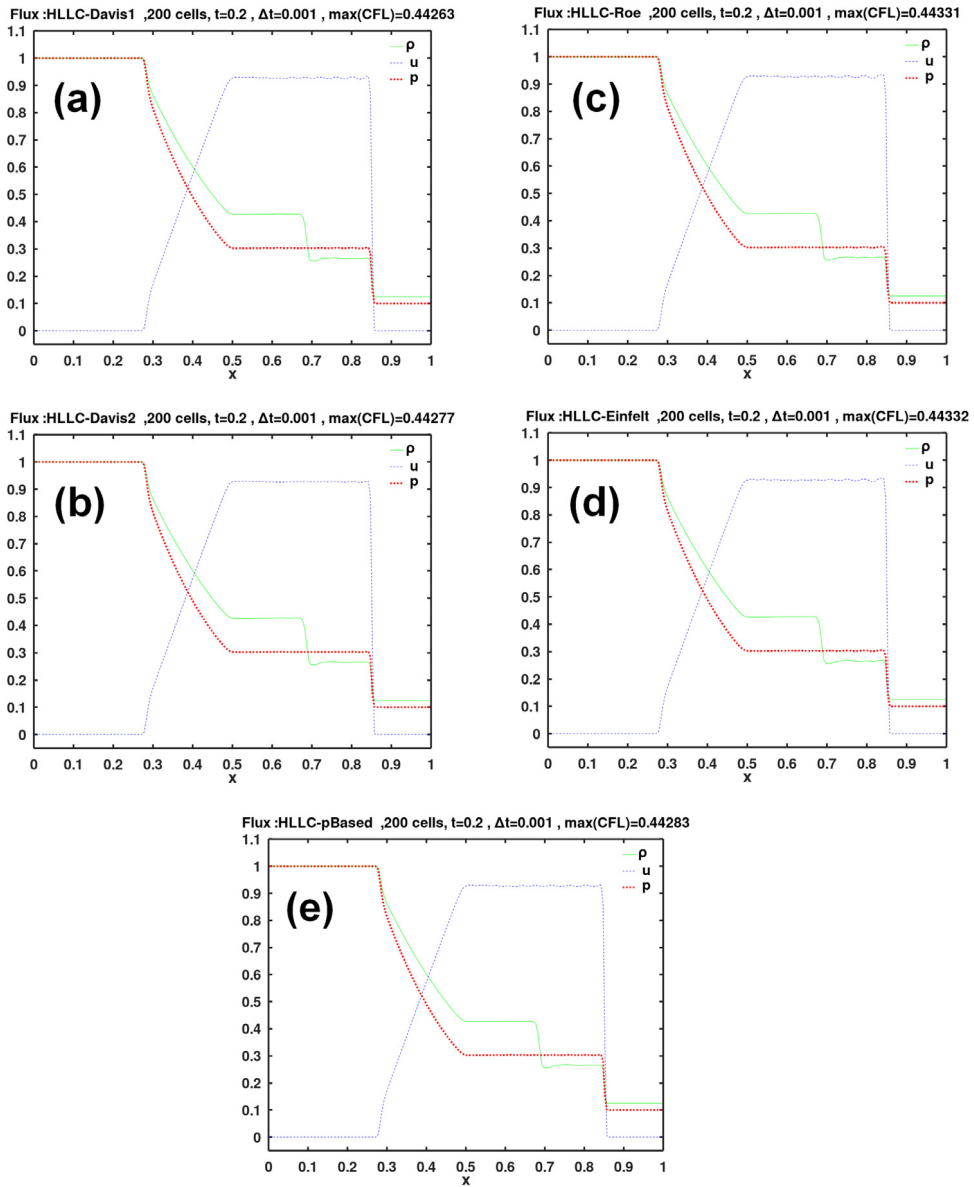

**Fig. 13.** Calculated final primitive variables using face flux of the Harten, Lax, van Leer, contact (HLLC) family of methods. The right/left wave speeds are obtained using: (a) first estimation of Davis, (b) second estimation of Davis, (c) eigenvalues of Reo, (d) eigenvalues of Einfelt, and (e) pressure-based calculation.





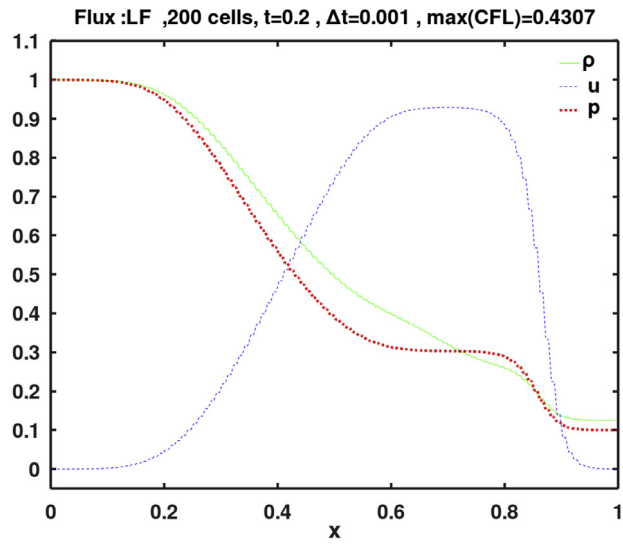

**Fig. 14.** Calculated final primitive variables using face flux of Lax and Friedrichs (LF).

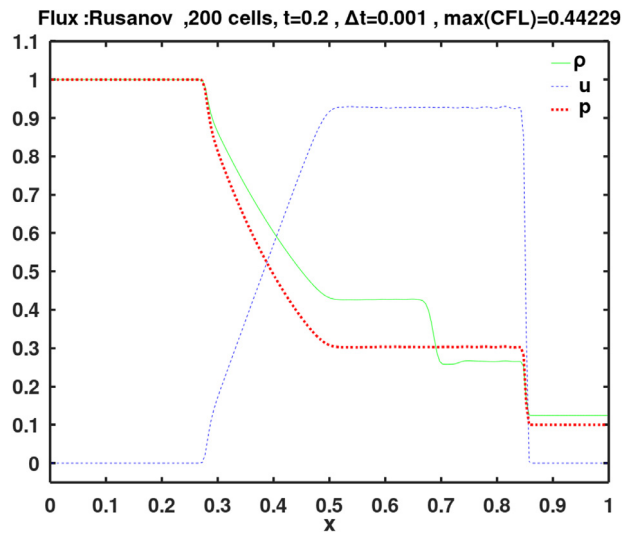

**Fig. 15.** Calculated final primitive variables using face flux of Rusanov.





**Table 4**

Comparison of the elapsed time (the runtime) for each of the flux method used here. Values are ordered from smallest (best) to largest (worst).

| Index | Flux method | Elapsed time (in seconds) | Percentage Increase of time (relative to the shortest one) |
|---|---|---|---|
| 5 | SW | 41.13235 | 0.0% |
| 7 | AUSM | 41.54838 | 1.0% |
| 3 | KNP | 41.68739 | 1.3% |
| 4 | KT | 41.69439 | 1.4% |
| 10 | AUFS | 41.90940 | 1.9% |
| 6 | vanLeer | 41.91540 | 1.9% |
| 17 | HLLC-Davis2 | 42.30342 | 2.8% |
| 11 | HLL-Davis1 | 42.42343 | 3.1% |
| 22 | Rusanov | 42.56944 | 3.5% |
| 8 | AUSM+ | 42.78445 | 4.0% |
| 12 | HLL-Davis2 | 42.85645 | 4.2% |
| 16 | HLLC-Davis1 | 42.97146 | 4.5% |
| 13 | HLL-Roe | 43.10346 | 4.8% |
| 15 | HLL-pBased | 43.16347 | 4.9% |
| 14 | HLL-Einfelt | 43.36448 | 5.4% |
| 18 | HLLC-Roe | 43.80651 | 6.5% |
| 20 | HLLC-pBased | 43.82451 | 6.5% |
| 9 | AUSM+-up | 44.00052 | 7.0% |
| 2 | Roe | 44.67956 | 8.6% |
| 19 | HLLC-Einfelt | 44.76256 | 8.8% |
| 21 | LF | 45.50260 | 10.6% |
| 1 | Riemann | 46.49466 | 13.0% |

## Declaration of Competing Interest



## Acknowledgment

Dr. Michael Zingale, Associate Professor at Stony Brook University in the Department of Physics and Astronomy (Stony Brook, New York, U.S.A.), is highly appreciated for the kind permission of using a figure of the Castro project for validation purpose here.

The author also deeply appreciates valuable support from the National Energy Technology Laboratory (U.S. Department of Energy).

## Funding

This work is not a funded project.

## Appendix A. One-dimensional Euler equations

The Euler equations governing a fluid motion neglects friction effects (viscosity) completely [43], leading to a reduced form of the more comprehensive Navier-Stokes equations. While the Euler formulation may appear narrowly restricted, there are many conditions where the effect of fluid viscosity can be safely neglected and this simplified formulation becomes accurate enough.

In a vector form, the one-dimensional Euler equation governing inviscid fluids is [44]

$$\frac{\partial q}{\partial t} + \frac{\partial F}{\partial x} = 0 \tag{A.1}$$

where ($q$) is a vector of the conserved variables, ($t$) is the time, ($F$) is a vector of flux (sum of convective flux in the x-axis and a pressure contribution). The conserved variables vector is

$$q = \left\{ \begin{array}{c} \rho \\ \rho u \\ \rho\, e_T \end{array} \right\} \tag{A.2}$$

where ($\rho$) is the fluid density, which is the mass per unit volume; ($u$) is the fluid velocity, so ($\rho\, u$) is the momentum per unit volume; ($e_T$) is the total specific internal energy, so ($\rho\, e_T$) is the total internal energy per unit volume. The density ($\rho$), velocity ($u$), and pressure ($p$) are referred to as the primitive variables.





The total specific internal energy is the specific internal energy ($e$) plus the specific kinetic energy (kinetic energy per unit mass), thus

$$e_T = e + \frac{1}{2}u^2 \tag{A.3}$$

The flux vector is

$$F = \left\{ \begin{array}{c} \rho \ u \\ \rho u^2 + p \\ (\rho \ e_T + p) \ u \end{array} \right\} = \left\{ \begin{array}{c} \rho \ u \\ \rho u^2 + p \\ \rho \ u \ h_T \end{array} \right\} \tag{A.4}$$

where ($h_T$) is the total specific enthalpy

$$h_T = \ e_T + \frac{p}{\rho} = e + \frac{1}{2}u^2 + \frac{p}{\rho} \tag{A.5}$$

The specific enthalpy is

$$h = e + \frac{p}{\rho} = h_T - \frac{1}{2}u^2 \tag{A.6}$$

For a thermally perfect gas (also known as an ideal gas), the pressure, density, and absolute temperature ($T$) are related by the following form of the ideal gas equation of state:

$$p = \rho \ R \ T \tag{A.7}$$

where ($R$) is the specific gas constant. Also, for a thermally perfect gas, the specific heat capacity at constant volume ($C_v$) and the specific heat capacity at constant pressure ($C_p$) both vary with the temperature only, and the following relation holds true:

$$C_p(T) = C_v(T) + R \tag{A.8}$$

The ratio of the specific heat capacities ($\gamma$), also known as an the adiabatic index, is

$$\gamma = \frac{C_p}{C_v} \tag{A.9}$$

The following equations apply to a thermally perfect gas

$$C_v(T) = \frac{R}{\gamma(T) - 1} \tag{A.10}$$

$$C_p(T) = \frac{\gamma(T) \ R}{\gamma(T) - 1} \tag{A.11}$$

$$a = \sqrt{\gamma \ \frac{p}{\rho}} = \sqrt{\gamma R \ T} \tag{A.12}$$

$$\frac{p}{\rho} = \frac{a^2}{\gamma} \tag{A.13}$$

$$h = e + \frac{1}{\gamma}a^2 \tag{A.14}$$

$$h_T = e_T + \frac{1}{\gamma}a^2 \tag{A.15}$$

$$\rho u^2 + p = \rho a^2 \left( M^2 + \frac{1}{\gamma} \right) \tag{A.16}$$

where ($a$) is the speed of sound in the gas (the sonic speed), and ($M = u/a$) is the Mach number. In this work, the Mach number is a signed scalar quantity. Like the scalar velocity ($u$), a positive value of $M$ indicates a flow motion in the positive x-axis (to the right), whereas a negative value of $M$ indicates a flow motion in the negative x-axis (to the left).

If a thermally perfect gas has also the specialty that its specific enthalpy and its specific internal energy can be related as

$$h(T) = \gamma \ e(T) \tag{A.17}$$

with ($\gamma$) being now a constant, independent of the temperature, then the gas is described as calorically perfect and its specific heat capacities are also independent of the temperature [45]. The calorically perfect gas assumption permits analytical analysis for subsonic (Mach number magnitude below unity) flows and supersonic (Mach number magnitude above unity) flows [46,47].





For a calorically perfect gas, the following formulas apply

$$p = (\gamma - 1)\rho\, e \tag{A.18}$$

$$e = \frac{p}{\rho\,(\gamma - 1)} = \frac{1}{\gamma\,(\gamma - 1)}a^2 \tag{A.19}$$

$$h = \frac{\gamma}{\gamma - 1}\,\frac{p}{\rho} = \frac{1}{\gamma - 1}a^2 \tag{A.20}$$

$$\gamma\, e_T = \frac{\gamma}{2}u^2 + \frac{1}{\gamma - 1}a^2 \tag{A.21}$$

$$h_T = \frac{1}{2}u^2 + \frac{1}{\gamma - 1}a^2 = \frac{1}{2}u^2 + \frac{\gamma}{\gamma - 1}\,\frac{p}{\rho} \tag{A.22}$$

$$\rho\, u\, h_T = \rho a^3\left(\frac{M^3}{2} + \frac{M}{\gamma - 1}\right) \tag{A.23}$$